\documentclass{article}
\usepackage{amsfonts}
\usepackage{amsthm}
\usepackage{mathrsfs}
\usepackage{amsmath}
\usepackage{url}

\theoremstyle{definition}
\newtheorem{theorem}{Theorem}
\newtheorem{lemma}[theorem]{Lemma}
\newtheorem{example}[theorem]{Example}
\newtheorem{proposition}[theorem]{Proposition}
\newtheorem{corollary}[theorem]{Corollary}
\newtheorem{definition}{Definition}

\newtheorem{remark}[theorem]{Remark}

\def\N{\mathbb{N}}

\def\Lbor{\mathscr{L}_{Bor}}
\def\Mf{\mathscr{M}_f}
\def\Apply{\mbox{Apply}}

\begin{document}

\title{The First-Order Syntax of Variadic Functions}
\author{Samuel Alexander\footnote{alexander@math.ohio-state.edu}}
\date{2013}

\maketitle








\begin{abstract}
 We extend first-order logic to include variadic function symbols, and 
 prove
 a substitution lemma.
 Two applications are given: one to bounded quantifier elimination and
 one to the definability of certain Borel sets.
\end{abstract}

\maketitle

\section{Introduction}

A variadic function is a function which takes a variable number of arguments:
for example, a function from $\N^{<\N}$ to $\N$ is variadic,
where $\N^{<\N}$ denotes the set of finite sequences of naturals.
In classical first-order logic, a language has function symbols of fixed arities.
In this paper I will explore how variadic function symbols can be added to
first-order logic.  In so doing, we will also formalize the syntax of the ellipsis,
$\cdots$, which of course is closely related to variadic function symbols.

To get an idea of the subtleties of the ellipsis, consider the following ``proof'' that $5050=385$:
\begin{enumerate}
\item We know $1+\cdots+100=5050$.
\item We know $1=1^2$ and $100=10^2$.
\item By replacement, $1^2+\cdots+10^2=5050$.
\item Also, $1^2+\cdots+10^2=(10)(10+1)(2\cdot 10+1)/6=385$. So $5050=385$.
\end{enumerate}

Evidently, mathematicians
implicitly impose some special syntax on the ellipsis.  This will be made explicit in
the paper.

Of course, we can already talk about unary functions $\N\to\N$ which interpret their input
as the code for a finite sequence.  My hope is that some coding can be avoided
by allowing variadic function symbols.

I was led to investigate the syntax of variadic function symbols when I was investigating
a certain class of subsets of Baire space and realized that I could characterize that class
with the help of first-order logic extended by variadic function symbols.  The results
are written up in \cite{alexander}.  Some of the basic results of this paper were first published
there.

Variadic functions are used in many programming languages.
What little literature presently exists mostly seems to be in this context (for example,
Byrd and Friedman \cite{byrd}) and in the related context of $\lambda$-calculus (for example Goldberg \cite{goldberg}).

\section{Basic definitions}
\label{sect2}

\begin{definition}
\hfill
\begin{itemize}
\item
A \emph{first-order variadic language} (or simply a \emph{variadic language}) is a first-order language,
including a constant symbol $\mathbf{n}$ for every $n\in\N$, together with a set of
\emph{variadic function symbols}, and a special symbol $\cdots_x$ for every variable $x$.
\item
A \emph{structure} (or a \emph{model}) for a variadic language $\mathscr{L}$ is a structure $\mathscr{M}$ for the first-order part
of $\mathscr{L}$, with universe $\N$,
and which interprets each $\mathbf{n}$ as $n$,
together with a set of variadic functions $\N^{<\N}\to\N$, one for every variadic function symbol of 
$\mathscr{L}$.
If $G$ is a variadic function symbol of $\mathscr{L}$, $G^{\mathscr{M}}$ will denote the corresponding function $G^{\mathscr{M}}:\N^{<\N}\to \N$.
\end{itemize}
\end{definition}

\begin{definition}
\hfill
\begin{itemize}
\item
\emph{(Terms)}
If $\mathscr{L}$ is a variadic language then the \emph{terms} of $\mathscr{L}$, and their free variables, are defined as follows:
\begin{enumerate}
\item
If $c$ is a constant symbol, then $c$ is a term with $FV(c)=\emptyset$.
\item
If $x$ is a variable, then $x$ is a term with $FV(x)=\{x\}$.
\item
If $f$ is an $n$-ary or variadic function symbol and $u_1,\ldots,u_n$ are terms,
then $f(u_1,\ldots,u_n)$ is a term with free variables $FV(u_1)\cup\cdots\cup FV(u_n)$.
\item
If $G$ is a variadic function symbol, $u,v$ are terms, and $x$ is a variable, then
\[
G(u(\mathbf{0})\cdots_x u(v))\]
is a term, with free variables $(FV(u)\backslash\{x\})\cup FV(v)$.
\end{enumerate}
\item
\emph{(Term substitution)}
If $r,t$ are terms and $x$ is a variable, then the term $r(x|t)$ obtained \emph{by substituting $t$ for $x$ in $r$} is
defined by induction in the usual way, with two new cases:
\begin{enumerate}
\item
If $r$ is $G(u(\mathbf{0})\cdots_x u(v))$ then $r(x|t)$ is
\[
G(u(\mathbf{0})\cdots_x u(v(x|t))).\]
\item
If $r$ is $G(u(\mathbf{0})\cdots_y u(v))$ where $y\not=x$, then $r(x|t)$ is
\[
G(u(x|t)(\mathbf{0})\cdots_y u(x|t)(v(x|t))).\]
\end{enumerate}
\end{itemize}
\end{definition}

\begin{lemma}
\label{uniquereadability}
\emph{(Unique Readability)}
Assume $\mathscr{L}$ has the following properties:
\begin{enumerate}
\item
Every symbol of $\mathscr{L}$ is \emph{exactly one} of the following:
a left parenthesis, a right parenthesis, a logical connective, =, a constant symbol, a variable, an $n$-ary predicate symbol for some $n$,
an $n$-ary function symbol for some $n$, a variadic function symbol, or an ellipsis $\cdots_x$ for some variable $x$.
\item
If some symbol is an $n$-ary function (resp.~predicate) symbol and also an $m$-ary function (resp.~predicate) symbol,
then $n=m$.
\item
If some symbol is $\cdots_x$ and $\cdots_y$, then $x$ and $y$ are the same variable.
\end{enumerate}
Then the terms of $\mathscr{L}$ have the unique readability property.
\end{lemma}

\begin{proof}
By the usual inductive argument.
\end{proof}

Henceforth, we will always assume every language satisfies the hypotheses of Lemma~\ref{uniquereadability}.
Thus the free variables of a term are well-defined, as is term substitution.

\begin{example}
\label{finitesummationnotation}
\emph{(Finite Sigma Notation)}  If $\mathscr{L}$ is a first-order language, we can extend it to a variadic language $\mathscr{L}_{\Sigma}$
by adding a variadic function symbol $\Sigma$ (along with ellipses and numerals).
In practice, the term $\Sigma(u(\mathbf{0})\cdots_x u(v))$ is often written:
\[
\sum_{x=0}^{v}u.\]
The obvious interpretation $\Sigma^{\mathscr{M}}$ of $\Sigma$ in a structure is the variadic addition function
$\Sigma^{\mathscr{M}}(a_0,\ldots,a_n)=a_0+\cdots+a_n$.
\end{example}

Throughout the paper, if $\mathscr{M}$ is a structure, then an \emph{assignment} shall mean a function which maps variables to elements
of the universe of $\mathscr{M}$.
If $s$ is an assignment
and $n\in\N$, I will write $s(x|n)$ for the assignment which agrees with $s$ everywhere except that it maps $x$
to $n$.

We would like to define the interpretation of a term in a model by an assignment.  This is straightforward in classic logic
but when variadic terms are introduced, interpretation becomes more subtle.
There are actually two possible definitions; they are equivalent, but to show it, we will first need to establish a substitution lemma
for one of the two.

When naively defining a numerical value for
$\sum_{i=0}^{v} u(i)$, where $u(i)$ and $v$ are mathematical expressions, we implicitly use a
definition by recursion on expression complexity, as each summand $u(i)$ may itself involve nested summations.  The process terminates
because each summand $u(i)$ is strictly simpler than $\sum_{i=0}^{10} u(i)$, which is true because $i$ itself
is not a compound expression but a natural number.  
Now, we'd like to say $\sum_{i=0}^{v} u(i) = u(0)+\cdots+u(v)$, where the summands on the right
are recursively computed using the definition currently being made.  But there are two ways to get here formally:
\begin{enumerate}
\item \emph{(Syntactic)} Write a list of $v+1$ terms, the $n$th of which is obtained by syntactically replacing the variable $i$ in $u$ by the constant 
$n$.
Recursively compute and add each of these new terms, using the same values for variables as we're currently using.
\item \emph{(Semantic)} Write a list of $v+1$ terms, each of which is exactly $u$.
Recursively compute and add them, but when computing the $n$th one, do it assuming the value of variable $i$ is $n$.
\end{enumerate}
This motivates the following definition of two interpretations (in Corollary~\ref{interpunambiguous} we will
see that both interpretations are equivalent).

\begin{definition}
\label{terminterp}
\emph{(Term interpretation)}
Let $\mathscr{M}$ be a structure for a variadic language $\mathscr{L}$, and let $s$ be an assignment.
Assume we've defined natural number interpretations $u^{s'},u_{s'}$ (respectively, \emph{syntactic} and \emph{semantic}
interpretations of $u$) for every assignment $s'$ and every term $u$ strictly simpler than $t$.
We define $t^s$ and $t_s$ inductively according to the following cases:
\begin{enumerate}
\item If $t$ is a constant symbol $c$, then $t^s=t_s=c^{\mathscr{M}}$.
\item If $t$ is a variable $x$, then $t^s=t_s=s(x)$.
\item If $t$ is $f(u_1,\ldots,u_k)$, then $t^s=f^{\mathscr{M}}(u_1^s,\ldots,u_k^s)$ and $t_s=f^{\mathscr{M}}({u_1}_s,\ldots,{u_k}_s)$.
\item If $t$ is $G(u(\mathbf{0})\cdots_x u(v))$, then
\begin{align*}
G(u(\mathbf{0})\cdots_x u(v))^s &= G^{\mathscr{M}}(u(x|\mathbf{0})^s,\ldots,u(x|\overline{v^s})^s),\\
G(u(\mathbf{0})\cdots_x u(v))_s &= G^{\mathscr{M}}(u_{s(x|0)},\ldots,u_{s(x|v_s)}).\end{align*}
Here $\overline{v^s}$ denotes the constant symbol corresponding to the natural $v^s$.
\end{enumerate}
\end{definition}

\begin{example}To illustrate the definition, assume $v^s=v_s=5$.
Then by definition
\begin{align*}
G(u(\mathbf{0})\cdots_x u(v))^s
&=
G^{\mathscr{M}}(u(x|\mathbf{0})^s,\ldots,u(x|\mathbf{5})^s)\\
&=G^{\mathscr{M}}(u(x|\mathbf{0})^s,u(x|\mathbf{1})^s,u(x|\mathbf{2})^s,u(x|\mathbf{3})^s,u(x|\mathbf{4})^s,u(x|\mathbf{5})^s).\\
G(u(\mathbf{0})\cdots_x u(v))_s
&=
G^{\mathscr{M}}(u_{s(x|0)},\ldots,u_{s(x|5)})\\
&=
G^{\mathscr{M}}(u_{s(x|0)},u_{s(x|1)},u_{s(x|2)},u_{s(x|3)},u_{s(x|4)},u_{s(x|5)}).
\end{align*}
\end{example}

\begin{remark}
Definition~\ref{terminterp} may seem somewhat suspect because
of how it uses meta-ellipses to define the semantics of ellipses.
If we were forbidden from using meta-ellipses to define the semantics
of ellipses, there are two approaches we could take.
One approach would be to use simultaneous induction to simultaneously
define interpretations $t^s$ and $t_s$ and also define sequences
$\langle t(x|\mathbf{0})^s,\ldots,t(x|\mathbf{n})^s\rangle$ and
$\langle t_{s(x|0)},\ldots,t_{s(x|n)}\rangle$.  The latter would be defined
by induction on $n$ by means of \emph{concatenation}:
\begin{align*}
\langle t(x|\mathbf{0})^s,\ldots,t(x|\overline{n+1})^s\rangle
&=
\langle t(x|\mathbf{0})^s,\ldots,t(x|\overline{n})^s\rangle \frown \langle t(x|\overline{n+1})^s\rangle\\
\langle t_{s(x|0)},\ldots,t_{s(x|n+1)}\rangle
&=
\langle t_{s(x|0)},\ldots,t_{s(x|n)}\rangle \frown \langle t_{s(x|n+1)}\rangle,\end{align*}
assuming that $(t')^{s'}$ and $(t')_{s'}$ are already defined for every assignment $s'$ and every term
$t'$ at most as complex as $t$; meanwhile, $t^s$ and $t_s$ would be defined as in Definition~\ref{terminterp}
except that we'd let
\begin{align*}
G(u(\mathbf{0})\cdots_x u(v))^s &= G^{\mathscr{M}}(\langle u(x|\mathbf{0})^s,\ldots,u(x|\overline{v^s})^s\rangle),\\
G(u(\mathbf{0})\cdots_x u(v))_s &= G^{\mathscr{M}}(\langle u_{s(x|0)},\ldots,u_{s(x|v_s)}\rangle).\end{align*}
This approach does not truly use meta-ellipses
except as a name; the name
could be changed without changing the definition.
Another alternative approach would be to use \emph{generalized structures} which we'll discuss in Section~\ref{sectmechanize}.
\end{remark}

\begin{remark}
The syntactic part of Definition~\ref{terminterp} relies on the fact that $u(x|\mathbf{c})$ is
strictly simpler than $G(u(\mathbf{0})\cdots_x u(v))$ for any constant symbol $\mathbf{c}$.
The minimalist might wonder whether we can treat
first-order variadic function symbols without so many constant symbols,
using only the constant symbol $\mathbf{0}$ and the successor function symbol $S$.
Maybe the natural way to translate the definition of $G(u(\mathbf{0})\cdots_x u(v))^s$
would be to take Definition~\ref{terminterp} using numerals $\mathbf{n}=SS\ldots S(\mathbf{0})$.
But then $u(x|\mathbf{c})$ would no longer necessarily be simpler than $G(u(\mathbf{0})\cdots_x u(v))$,
casting doubt on the productiveness of the definition.
One way around this dilemma would be to define term complexity not as a natural number but
as an ordinal in $\epsilon_0$, defining the complexity of $G(u(\mathbf{0})\cdots_x u(v))$
to be, say, $\omega^{c(u)+c(v)}$, where $c(u)$ and $c(v)$ are the complexities of $u$ and $v$.
Of course, such a radical approach is not \emph{necessary}, but it is more elegant than other
solutions to the dilemma, and this author considers it a nice and unexpected application of ordinals to
syntax.
\end{remark}

\section{The Substitution Lemma}
\label{sect3}

We will deal mainly with syntactic interpretations $t^s$.  We will obtain a substitution lemma
for these, and use it to show that the two interpretations are identical.  The choice is arbitrary:
one could also obtain a substitution lemma about semantic interpretations and use that to show equality.
Once either version of the substitution lemma is obtained, and the two interpretations are shown equal,
the other substitution lemma becomes trivial.  In any case, technical lemmas are required.

\begin{lemma}
\label{techlemma1}
Suppose $u,t$ are terms and $x,y$ are variables.
\begin{enumerate}
\item
If $x$ is not a free variable of $u$, then $u(x|t)=u$.
\item
If $x$ is not a free variable of $u$, and $s$ is an assignment, then $u^s$ does not depend on $s(x)$.
\item
If $x$ is not a free variable of $u$ or $t$, then $x$ is not a free variable of $u(y|t)$.
\end{enumerate}
\end{lemma}

\begin{proof}
A straightforward induction.
\end{proof}

In first-order logic, substitutability is a property of formulas, but it is not needed for terms: if $r,t$ are any terms and $x$ is any variable,
then $t$ is substitutable for $x$ in $r$ (in first-order logic).  This breaks down in variadic logic, requiring a notion of substitutability into
terms.

\begin{definition}
\label{termsubstitutable}
Fix a term $t$ and a variable $x$.  The substitutability of $t$ for $x$ in a term $r$ is defined inductively:
\begin{itemize}
\item
If $r$ is a variable or constant symbol, then $t$ is substitutable for $x$ in $r$.
\item
If $r$ is $f(u_1,\ldots,u_n)$ where $u_1,\ldots,u_n$ are terms and $f$ is an $n$-ary or variadic function symbol,
then $t$ is substitutable for $x$ in $r$ if and only if $t$ is substitutable for $x$ in all the $u_i$.
\item
If $r$ is $G(u(\mathbf{0})\cdots_y u(v))$ where $G$ is a variadic function symbol, $u,v$ are terms, and $y$ is a variable
(which may or may not be $x$), then $t$ is substitutable for $x$ in $r$
if at least one of the following hold:
\begin{itemize}
\item
$x$ is not a free variable of $r$, or
\item
$y=x$ and $t$ is substitutable for $x$ in $v$, or
\item
$y$ is not a free variable of $t$ and $t$ is substitutable for $x$ in both $u$ and $v$.
\end{itemize}
\end{itemize}
\end{definition}

For a non-substitutability example, consider the term $\sum_{y=0}^{\mathbf{10}}x\cdot y$ and try
substituting $t=y$ for $x$.  The result is $\sum_{y=0}^{\mathbf{10}}y\cdot y$, which is no good
since the new occurrence of $y$ becomes bound by the summation.

We define substitutability for formulas in the usual way, with just one change:
if $p$ is a predicate symbol (or $=$), and $u_1,\ldots,u_n$ are terms,
we say $t$ is substitutable for $x$ in $p(u_1,\ldots,u_n)$ if and only if $t$ is substitutable for $x$ in each $u_i$.

\begin{lemma}
\label{techlemma1point5}
Suppose $t$ is a term which is substitutable for the variable $x$ in $G(u(\mathbf{0})\cdots_y u(v))$.
Then $t$ is substitutable for $x$ in $v$.
And
if $x\not=y$, then $t$ is substitutable for $x$ in $u$.
\end{lemma}

\begin{proof}By induction.
\end{proof}

\begin{lemma}
\label{techlemma2}
Suppose $r,t$ are terms, $x,y$ are distinct variables, and $c$ is a constant symbol.
If $t$ is substitutable for $x$ in $r$, and $y$ does not occur free in $t$, then $r(x|t)(y|c)=r(y|c)(x|t)$.
\end{lemma}

\begin{proof}
By induction on complexity of $r$.  If $r$ is a constant symbol or $f(r_1,\ldots,r_n)$
for some function symbol $f$ and terms $r_1,\ldots,r_n$ (wherein $t$ is substitutable for $x$),
the lemma is clear by induction.  If $r$ is a variable,
the claim follows since $y$ does not
occur free in $t$.
But suppose $r$ is $G(u(\mathbf{0})\cdots_z u(v))$ for some variadic function symbol $G$,
terms $u,v$,
and variable $z$ (which may be $x$, $y$, or neither).
Since $t$ is substitutable for $x$ in $r$, at least one of the following holds:
$x$ is not free in $r$; or $z=x$ and $t$ is substitutable for $x$ in $v$; or $t$ is substitutable for $x$ in $u$ and $v$.
If $x$ is not free in $r$, then the lemma follows from Lemma~\ref{techlemma1}.
But suppose $x$ is free in $r$.
Unravelling definitions:

\begin{figure}[h]
\centering
\begin{tabular}{|l|l|}
\hline
If $z$... & Then $r(x|t)$ equals...\\
\hline
\hline
$=x$
&
$G(u(\mathbf{0})\cdots_z u(v(x|t)))$\\
\hline
$\not=x$
&
$G(u(x|t)(\mathbf{0})\cdots_z u(x|t)(v(x|t)))$\\
\hline
\hline
If $z$... & Then $r(x|t)(y|c)$ equals...\\
\hline
\hline
$=x$
&
$G(u(y|c)(\mathbf{0})\cdots_z u(y|c)(v(x|t)(y|c)))$\\
\hline
$=y$
&
$G(u(x|t)(\mathbf{0})\cdots_z u(x|t)(v(x|t)(y|c)))$\\
\hline
$\not\in\{x,y\}$
&
$G(u(x|t)(y|c)(\mathbf{0})\cdots_z u(x|t)(y|c)(v(x|t)(y|c)))$\\
\hline
\end{tabular}
\end{figure}

By Lemma~\ref{techlemma1point5}, $t$ is substitutable for $x$ in $v$, so $v(x|t)(y|c)=v(y|c)(x|t)$
by induction.
And if $z\not=x$, then Lemma~\ref{techlemma1point5} tells us $t$ is substitutable for $x$ in $u$ as well,
and so by induction $u(x|t)(y|c)=u(y|c)(x|t)$.
The lemma follows by using these facts to rewrite the last row of the table
and compare with a similar table for $r(y|c)(x|t)$.
\end{proof}

\begin{theorem}
\label{substitutionlemma}
\emph{(The Variadic Substitution Lemma for Terms)}
Let $\mathscr{M}$ be a structure for a variadic language $\mathscr{L}$ and let $s$
be an assignment.
If $r$ and $t$ are terms such that $t$ is substitutable for $x$ in $r$,
then $r(x|t)^s=r^{s(x|t^s)}$.
\end{theorem}

\begin{proof}
We induct on the complexity of $r$, and most cases
are straightforward.  If $x$ is not free in $r$, the claim is trivial; assume
$x$ is free in $r$.  The two important cases follow.

We must show $G(u(\mathbf{0})\cdots_y u(v))(x|t)^s=G(u(\mathbf{0})\cdots_y u(v))^{s(x|t^s)}$ when $y$ is a 
different variable than $x$ and $t$ is substitutable for $x$ in $G(u(\mathbf{0}) \cdots_y u(v))(x|t)$.
Using induction:
\begin{align*}
G(u(\mathbf{0})\cdots_y u(v))(x|t)^s
&=
G(u(x|t)(\mathbf{0}) \cdots_y u(x|t)(v(x|t)))^s\\
&=
G^{\mathscr{M}}\left(u(x|t)(y|\mathbf{0})^s,\ldots,u(x|t)\left(y\left|\overline{v(x|t)^s}\right.\right)^s\right)
\\
&=
G^{\mathscr{M}}\left(u(y|\mathbf{0})(x|t)^s,\ldots,u\left(y\left|\overline{v(x|t)^s}\right.\right)(x|t)^s\right)
& \mbox{($*$)}
\\
&=
G^{\mathscr{M}}
\left(u(y|\mathbf{0})^{s(x|t^s)},\ldots,
u\left(y
\left|
\overline{v^{s(x|t^s)}}
\right.\right)^{s(x|t^s)}
\right)
\\
&=
G(u(\mathbf{0}) \cdots_y u(v))^{s(x|t^s)}.
\end{align*}

To reach line $(*)$, we need the fact that $u(x|t)(y|c)=u(y|c)(x|t)$ for any constant symbol $c$.
If $x$ is not free in $r$, then it's not free in $u$ (since $y\not=x$),
and so this follows from Lemma~\ref{techlemma1}.
Otherwise, since $t$ is substitutable for $x$ in $r$ by Lemma~\ref{techlemma1point5},
we must have that $y$ does not occur free in $t$, and so we can invoke Lemma~\ref{techlemma2}.

We must also show 
$G(u(\mathbf{0})\cdots_x u(v))(x|t)^s=G(u(\mathbf{0}) \cdots_x u(v))^{s(x|t^s)}$ when 
$t$ is substitutable for $x$ in $G(u(\mathbf{0}) \cdots_x u(v))(x|t)$.
Using induction:
\begin{align*}
G(u(\mathbf{0}) \cdots_x u(v))(x|t)^s
&=
G(u(\mathbf{0}) \cdots_x u(v(x|t)))^s\\
&=
G^{\mathscr{M}}\left(
u(x|\mathbf{0})^s
,\ldots,u
\left(x
\left|\overline{
v(x|t)^s
}
\right.\right)^s\right)\\
&=
G^{\mathscr{M}}\left(u(x|\mathbf{0})^s,\ldots,u\left(x\left|\overline{v^{s(x|t^s)}}\right.\right)^s\right)\\
&=
G^{\mathscr{M}}\left(u(x|\mathbf{0})^{s(x|t^s)},\ldots,u\left(x\left|\overline{v^{s(x|t^s)}}\right.\right)
^{s(x|t^s)}\right)
&
\mbox{($**$)}\\
&=
G(u(\mathbf{0}) \cdots_x u(v))^{s(x|t^s)}.
\end{align*}
In line $(**)$, I am able to change ``exponents'' from $s$ to $s(x|t^s)$ because
the terms in question do not depend on $x$.
\end{proof}

\begin{corollary}
\label{interpunambiguous}
For any term $t$ and assignment $s$, $t^s=t_s$.
\end{corollary}

\begin{proof}
By induction on $t$.  All cases are immediate except the case when $t$ is $G(u(\mathbf{0})\cdots_x u(v))$.
Note that constant symbols are always substitutable and write:
\begin{align*}
G(u(\mathbf{0})\cdots_x u(v))^s
&= G^{\mathscr{M}}\left(u(x|\mathbf{0})^s,\ldots,u(x|\overline{v^s})^s\right)
& \mbox{(By Definition)}\\
&=
G^{\mathscr{M}}\left(u^{s(x|0)},\ldots,u^{s(x|v^s)}\right)
& \mbox{(By Theorem~\ref{substitutionlemma})}\\
&=
G^{\mathscr{M}}\left(u_{s(x|0)},\ldots,u_{s(x|v_s)}\right)
& \mbox{(By Induction)}\\
&=
G(u(\mathbf{0})\cdots_x u(v))_s. & \mbox{(By Definition)}\end{align*}
\end{proof}

First-order formulas over a variadic language are now defined in the obvious way.
By Corollary~\ref{interpunambiguous}, we can 
define $\mathscr{M}\models t=r[s]$ if and only if $t^s=r^s$, or equivalently $t_s=r_s$;
that is, we are saved from having to make an arbitrary choice.
The remaining semantics are defined inductively in exactly the
same way they are for first-order logic.
If $\mathscr{M}$ is a structure, $s$ an assignment, and $\phi$ a formula,
then $\mathscr{M}\models\phi[s]$ is defined in the usual way from the above
atomic case, and $\mathscr{M}\models\phi$ means that $\mathscr{M}\models\phi[s]$
for every assignment $s$.
Term substitution in a formula is defined as usual.
Substitutability of a term for a variable in a formula is defined as usual,
except that in the atomic case, we say $t$ is substitutable for $x$ in $r=q$
if and only if $t$ is substitutable for $x$ in $r$ and $q$ (in the sense of Definition~\ref{termsubstitutable}).

\begin{corollary}
\emph{(The Variadic Substitution Lemma)}
If $t$ is a term which is substitutable for the variable $x$ in the formula $\phi$,
and $s$ is an assignment and $\mathscr{M}$ a structure, then $\mathscr{M}\models \phi(x|t)[s]$
if and only if $\mathscr{M}\models \phi[s(x|t^s)]$.
\end{corollary}

\begin{proof}
The proof is identical to the proof of the first-order substitution lemma, except that
Theorem~\ref{substitutionlemma} is invoked for the atomic case.
\end{proof}




\begin{example}
Working in an appropriate language and structure, it can be shown that
\[
\sum_{x=0}^{x}x = x(x+\mathbf{1})/\mathbf{2},\]
showing that it is safe to use the same variable in different roles, so long
as we use Definition~\ref{terminterp} to be completely clear what the
truth of the formula means.
\end{example}

\section{Bounded Quantifier Elimination}
\label{sect4}

In this section, we shall assume our languages have no predicate symbols.
If a language has a binary function symbol $\leq$ and a constant symbol $\mathbf{1}$, I will write $u\leq v$
to abbreviate $\leq(u,v)=\mathbf{1}$.

\begin{definition}
\hfill
\begin{itemize}
\item
If $\mathscr{L}$ is a variadic language, the \emph{quantifier-free} formulas of $\mathscr{L}$ are defined inductively:
$\phi$ is quantifier-free whenever $\phi$ is atomic; and if $\phi$ and $\psi$ are quantifier-free, then so are
$\phi\wedge\psi$, $\phi\vee\psi$, $\phi\rightarrow\psi$, $\phi\leftrightarrow\psi$, and $\neg\phi$.
\item
If $\mathscr{L}$ contains a binary function symbol $\leq$ and constant symbol $\mathbf{1}$,
the \emph{unbounded-quantifier-free} (or uqf) formulas of $\mathscr{L}$ are defined 
inductively:
$\phi$ is uqf whenever $\phi$ is atomic; and if $\phi$ and $\psi$ are uqf and $x,y$ are distinct variables,
then $\phi\wedge\psi$, $\phi\vee\psi$, $\phi\rightarrow\psi$, $\phi\leftrightarrow\psi$, $\neg\phi$,
$\exists x\,(x\leq y\wedge \phi)$, and $\forall x\,(x\leq y\rightarrow\phi)$ are also uqf.
\end{itemize}
\end{definition}

\begin{proposition}
\label{qelim}
\emph{(Bounded Quantifier Elimination)}
Suppose $\mathscr{L}$ is a variadic language containing (possibly among other things) binary
function symbols $\leq$, $+$ and $\delta$, and
a variadic function symbol $G$.
Suppose $\mathscr{M}$ is an $\mathscr{L}$-model which interprets $+$ as addition and interprets $\leq$, $\delta$, and $G$ by
\begin{align*}
\leq^{\mathscr{M}}(m,n) &= \left\{\begin{array}{l} \mbox{$1$ if $m\leq n$,}\\ \mbox{$0$ otherwise}\end{array}\right.
&
\delta^{\mathscr{M}}(m,n) &= \left\{\begin{array}{l} \mbox{$1$ if $m=n$,}\\ \mbox{$0$ otherwise}\end{array}\right.\end{align*}
\begin{align*}
G^{\mathscr{M}}(m_0,\ldots,m_n) &= \left\{\begin{array}{l}\mbox{$1$ if $m_i\not=0$ for some $0\leq i\leq n$,}\\
\mbox{$0$ otherwise.}\end{array}\right.
\end{align*}
For any uqf $\mathscr{L}$-formula $\phi$, there is a quantifier-free $\mathscr{L}$-formula $\psi$,
with the same free variables as $\phi$,
such that $\mathscr{M}\models\phi\leftrightarrow\psi$.
\end{proposition}

\begin{proof}
I will show more strongly that for any uqf formula $\phi$, there is a term $t_{\phi}$,
with exactly the free variables of $\phi$,
such that $\mathscr{M}\models \phi\leftrightarrow (t_{\phi}=\mathbf{1})$
and $\mathscr{M}\models\neg\phi\leftrightarrow (t_{\phi}=\mathbf{0})$.
This is by induction on $\phi$:
\begin{itemize}
\item
If $\phi$ is $u=v$, take $t_{\phi}=\delta(u,v)$.
\item
If $\phi$ is $\psi\wedge\rho$, take $t_{\phi}=\delta(t_{\psi}+t_{\rho},\mathbf{2})$.
\item
If $\phi$ is $\neg\psi$, take $t_{\phi}=\delta(t_{\psi},\mathbf{0})$.
\item
If $\phi$ is $\exists x\,(x\leq y\wedge \psi)$ where $x\not=y$, take $t_{\phi}=G(t_{\psi}(\mathbf{0})\cdots_x t_{\psi}(y))$.
\item All other cases for $\phi$ are reduced to the above by basic logic (there is no predicate case by assumption).
\end{itemize}
In all but the $\exists$ case, it is routine to check $\mathscr{M}\models\phi\leftrightarrow(t_{\phi}=\mathbf{1})$,
$\mathscr{M}\models\neg\phi\leftrightarrow(t_{\phi}=\mathbf{0})$.
The $\exists$ case goes as follows.
Assume $\phi$ is $\exists x\,(x\leq y\wedge\psi)$, $y\not=x$.
By induction, $\mathscr{M}\models \psi\leftrightarrow (t_{\psi}=\mathbf{1})$
and $\mathscr{M}\models\neg\psi\leftrightarrow (t_{\psi}=\mathbf{0})$.
So $\mathscr{M}\models\psi\leftrightarrow (t_{\psi}=\mathbf{1})[s]$ for every assignment $s$.
Let $s$ be an assignment.  Then:
\begin{align*}
\mathscr{M} &\models \mbox{$\phi\,[s]$ iff}\\
\mathscr{M} &\models \mbox{$\exists x\, (x\leq y\wedge\psi)\,[s]$ iff}\\
\mbox{$\exists n\in\N$ s.t.~$\mathscr{M}$} &\models \mbox{$x\leq y \wedge \psi\,[s(x|n)]$ iff}\\
\mbox{$\exists n\leq s(y)$ s.t.~$\mathscr{M}$} &\models \mbox{$t_{\psi}=\mathbf{1}\,[s(x|n)]$ iff} & \mbox{($*$)}\\
\mbox{$\exists n\leq s(y)$ s.t.} &\hphantom{\models} \mbox{$t_{\psi}^{s(x|n)}=1$ iff}\\
G^{\mathscr{M}}\left(t_{\psi}^{s(x|0)},\ldots,t_{\psi}^{s(x|s(y))}\right) &= \mbox{$1$ iff}\\
G^{\mathscr{M}}\left(t_{\psi}(x|\mathbf{0})^s,\ldots,t_{\psi}\left(x|\overline{y^s}\right)^s\right) &= \mbox{$1$ iff} & \mbox{($**$)}\\
\mathscr{M} &\models G(t_{\psi}(\mathbf{0})\cdots_x t_{\psi}(y))=\mathbf{1}\,[s].
\end{align*}
In line $(*)$ we use the fact $s(x|n)(y)=s(y)$ since $y\not=x$.
In line $(**)$ we invoke the Variadic Substitution Lemma (noting that constant symbols are always substitutable for $x$).
\end{proof}

\begin{corollary}
Let $\mathscr{L}$ be the language with constant symbols $\mathbf{n}$ for all $n\in\N$,
binary function symbols $+$, $\cdot$, $\delta$ and $\leq$, and a variadic function symbol $G$.
Let $\mathscr{M}$ be the model which interprets
everything in the obvious way (interpreting $\delta$ and $G$ as above).
A set $X\subseteq\N$ is computably enumerable if and only if there is a quantifier-free
$\mathscr{L}$-formula $\phi$, with free variables a subset of $\{x,y\}$,
such that for all $n\in\N$, $n\in X\leftrightarrow \mathscr{M}\models \exists y\,\phi(x|\mathbf{n})$.
\end{corollary}

\begin{proof}
Let $\mathscr{L}_0=\mathscr{L}\backslash\{G\}$ be the first-order part of $\mathscr{L}$.
Assume $X\subseteq\N$ is c.e.
By computability theory, there is a uqf formula $\phi_0$ of 
$\mathscr{L}_0$
with the desired properties.
The corollary follows by Bounded Quantifier Elimination.
The converse is clear by Church's Thesis.
\end{proof}

\section{Defining Borel Sets}

I will further extend the concept of variadic function symbols, and apply the idea to show that
a certain language can define any $\mathbf{\Sigma}_n^0$ or $\mathbf{\Pi}_n^0$ subset of $\N^{\N}$ with a formula
of complexity $\Sigma_n$ or $\Pi_n$ (respectively), in a rather nice way.
My interest in using powerful language to define Borel sets is partially influenced by Vanden Boom \cite{vandenboom}, pp.~276-277.

By an \emph{extended variadic language} I mean a first-order language together with
various \emph{$n$-ary-by-variadic} function symbols (for various $n\geq 0$),
as well as constant symbols $\mathbf{n}$ for every $n\in\N$ and ellipses $\cdots_x$.
A structure for an extended variadic language is a structure $\mathscr{M}$ for the first-order part,
with universe $\N$ and which interprets each $\mathbf{n}$ as $n$,
together with a function $G^{\mathscr{M}}:\N^n\times\N^{<\N}\to\N$ for every $n$-ary-by-variadic function
symbol $G$.
Terms, term substitution, term interpretation, and term substitutability are defined in ways very similar
to our work in Section 2, and the Variadic Substitution Lemma is proved in almost an identical way.

\begin{definition}
Let $\langle \rangle$ be the empty sequence.
\begin{itemize}
\item
By $\Lbor$ I mean the extended variadic language with a special unary function symbol $\mathbf{f}$ along with, for every $n>0$
and every $\iota:\N^n\to\N^{<\N}\backslash\{\langle\rangle\}$,
an $n$-ary function symbol $\ell_{\iota}$ and an $n$-ary-by-variadic function symbol
$\tau_{\iota}$.
\item
For any $f:\N\to\N$, $\Mf$ is the $\Lbor$ structure which interprets $\mathbf{f}$ as $f$ and
which, for any $n>0$ and $\iota:\N^n\to\N^{<\N}\backslash\{\langle\rangle\}$, interprets
\begin{align*}
\ell_{\iota}^{\Mf}(a_1,\ldots,a_n) &= \mbox{the length of $\iota(a_1,\ldots,a_n)$, minus $1$}\\
\tau_{\iota}^{\Mf}(a_1,\ldots,a_n,b_1,\ldots,b_m) &=
\left\{\begin{array}{l}
\mbox{$1$ if $(b_1,\ldots,b_m)=\iota(a_1,\ldots,a_n)$,}\\
\mbox{$0$ otherwise.}\end{array}\right.
\end{align*}
\item
If $\phi$ is an $\Lbor$-sentence and $S\subseteq\N^{\N}$, say that $\phi$ \emph{defines} $S$ if,
for every $f:\N\to\N$, $f\in S$ if and only if $\Mf\models\phi$.
\end{itemize}
\end{definition}

\begin{theorem}
Let $n>0$ and let $S\subseteq\N^{\N}$.
Then $S$ is $\mathbf{\Sigma}_n^0$ (resp.~$\mathbf{\Pi}_n^0$) if and only if $S$ is defined by a $\Sigma_n$ (resp.~$\Pi_n$) sentence of $\Lbor$.
\end{theorem}

\begin{proof}
Obvious if $S=\emptyset$ or $S=\N^{\N}$, assume not.
If $f_0$ is a finite sequence of naturals, I'll write $[f_0]$ for the set of infinite extensions of $f_0$.

($\Rightarrow$) Assume $S$ is $\mathbf{\Sigma}_n^0$.
If $n$ is odd, we can write $S=\cup_{i_1\in\N}\cdots\cup_{i_n\in\N} [f_{i_1\cdots i_n}]$ where the $\{f_{i_1\cdots i_n}\}$ are
finite, nonempty sequences.  If $n$ is even, we can write $S=\cup_{i_1\in\N}\cdots\cap_{i_n\in\N} [f_{i_1\cdots i_n}]^c$.
Let $\iota:\N^n\to\N^{<\N}\backslash\{\langle\rangle\}$ be the map which sends $(i_1,\ldots,i_n)$ to $f_{i_1\cdots i_n}$.
Let $f:\N\to\N$.  For any $(i_1,\ldots,i_n)$, $f$ extends $f_{i_1\cdots i_n}$ if and only
if $\tau_{\iota}^{\Mf}(i_1,\ldots,i_n,f(0),\ldots,f(\ell_{\iota}^{\Mf}(i_1,\ldots,i_n)))=1$.
So if $n$ is odd, then $f\in S$ iff
\[
\Mf\models \exists x_1\cdots \exists x_n \tau_{\iota}(x_1,\ldots,x_n,\mathbf{f}(z)(\mathbf{0})\cdots_z 
\mathbf{f}(z)(\ell_{\iota}(x_1,\ldots,x_n)))=\mathbf{1}.\]
And if $n$ is even, then $f\in S$ iff
\[
\Mf\models \exists x_1\cdots \forall x_n \tau_{\iota}(x_1,\ldots,x_n,\mathbf{f}(z)(\mathbf{0})\cdots_z 
\mathbf{f}(z)(\ell_{\iota}(x_1,\ldots,x_n)))=\mathbf{0}.\]
The $\mathbf{\Pi}_n^0$ case is similar.

($\Leftarrow$) By induction on $n$.
For the base case, first use an induction argument on formula complexity to show that if $\phi$
is a quantifier-free sentence of $\Lbor$ and $\Mf\models\phi$ then there is some $k$ so big that
whenever $g:\N\to\N$ extends $(f(0),\ldots,f(k))$, then $\mathscr{M}_g\models\phi$.
Thus a set defined by a quantifier-free formula is open, hence clopen since its complement is also
defined by that formula's negation.
The base case follows: for example, if $S$ is defined by a sentence $\exists x\,\phi$,
then
(by Variadic Substitution) $S=\cup_{i\in\N}\{g:\N\to\N\,:\,\mathscr{M}_g\models \phi(x|\mathbf{i})\}$,
a $\mathbf{\Sigma}_1^0$ set since each unionand is clopen.
The induction case is straightforward.
\end{proof}

\section{A Partial Mechanization}
\label{sectmechanize}

We partially automated Sections~\ref{sect2} and \ref{sect3} using the
Coq proof assistant \cite{coq}.
In Coq, it is easier to work with the \emph{semantic}, rather than the \emph{syntactic},
term interpretations
of Definition~\ref{terminterp}.
This is because semantic term interpretation is recursive in a direct structural
way: to interpret a term, one needs only interpret direct subterms.
This is in contrast with syntactic term interpretation, which is recursive in term
depth.  To syntactically interpret a term, one must interpret terms which are not direct
subterms.  This makes it much more tedious to automate proofs about syntactic interpretations,
so our automation primarily deals with semantic interpretations.
We do, however,
automate Corollary~\ref{interpunambiguous}, in light of which, the distinction disappears.

Very often when automating mathematics, it is actually easier to prove a stronger result.
This is certainly the case here.  By a \emph{generalized structure} $\mathscr{M}$ for a variadic language
$\mathscr{L}$ we mean a structure for the first-order part of $\mathscr{L}$,
together with a set of interpretations $G^{\mathscr{M}}:\N^{\N}\times\N\to\N$ for each variadic function symbol
$G$ of $\mathscr{L}$.
This is a generalization in an obvious way: given a structure $\mathscr{M}'$ as in Section~\ref{sect2},
there corresponds a generalized structure $\mathscr{M}$ which agrees with $\mathscr{M}'$ on the first-order
part of $\mathscr{L}$ and is otherwise defined by
\[
G^{\mathscr{M}}(f,v) = G^{\mathscr{M}'}(f(0),\ldots,f(v))\]
whenever $G$ is a variadic function symbol, $v\in\N$, and $f\in\N^{\N}$.
The syntactic and semantic interpretations in $\mathscr{M}$
of a term $G(u(\mathbf{0})\cdots_x u(v))$
by an assignment $s$ are, respectively,
\begin{align*}
G(u(\mathbf{0})\cdots_x u(v))^s &= G^{\mathscr{M}}(k\mapsto u(x|\overline{k})^s,v^s),\\
G(u(\mathbf{0})\cdots_x u(v))_s &= G^{\mathscr{M}}(k\mapsto u_{s(x|k)},v_s).
\end{align*}
All the results of Sections~\ref{sect2} and \ref{sect3} generalize accordingly.
It is easier to automate these stronger results because Coq has better built-in
support for working with functions $\N\to\N$ than for working with finite sequences.

For even further simplicity, we also assume that all functions are either variadic or binary,
we assume the special constant symbols $\overline{c}$ are the only constant symbols in the
language, and we assume there are no predicate symbols.

Syntactic term interpretation seems to lie on the border of what Coq can handle.
Coq cannot automatically detect that the definition is total.  We are able to convince Coq
of its totality using an experimental feature of Coq called Program Fixpoint (Sozeau \cite{sozeau}).
Chung-Kil Hur (\cite{hur2011a}, \cite{hur2011b}) helped us
tremendously with the details of getting Program Fixpoint to work.

In performing this partial mechanization, we were influenced by R.~O'Connor's mechanization
of ordinary first-order logic \cite{oconnor}.

\section{Future Work}

There are several directions to take this study from here.
For one thing, Sections~\ref{sect2} and \ref{sect3} could
easily be extended to other types of logic.  In order to inject variadic
terms into a logic, there are two basic requirements: first, that function
terms make any sense at all in that logic; second, that the logic has a
semantics which plays well with variadic function symbols, especially the
ellipsis.  Some potential logics where we could add variadic function symbols
include second-order logic, more general multi-sorted logic, and nominal
logic, just to name three.  The question is not so much whether the machinery
can be added to the logic, but rather, what interesting applications result?

In the direction of multi-sorted logic, we could
deal with semantics where one sort ranges over (say) $\mathbb{R}$ and another
ranges over $\N$, and thereby rigorously study variadic functions living in the
real numbers (single-sorted first-order logic falls short here: how are we to
interpret a term like $\sum_{i=0}^{\pi}i$?)

One of the shortcomings of this first-order treatment is that we were not
able to give what should be a basic example: the general $\Apply$ function from
computer science.  If $G:\N^{<\N}\to\N$ is a variadic function and $n_1,\ldots,n_k\in\N$
then $\Apply(G,n_1,\ldots,n_k)$ is defined to be $G(n_1,\ldots,n_k)$.
This could be formalized using our variadic machinery in various typed logics
where it makes sense to have a function symbol whose ``arity'' is some cartesian product
of types.

Another direction we can go from here is to consider
function symbols of \emph{infinite arity}.  The basic idea is that if $G$
is an infinitary function symbol in a language and $u$ is a term and $x$ a variable,
then $G(u(\mathbf{0})\cdots_x)$ is another term, whose intended interpretation by
a model $\mathscr{M}$ and assignment $s$ is
\[
G(u(\mathbf{0})\cdots_x)^s = G^{\mathscr{M}}(u(x|\mathbf{0})^s,u(x|\mathbf{1})^s,\ldots),\]
where $G$ itself is interpreted as some infinitary $G^{\mathscr{M}}:\N^{\N}\to\N$.
In fact, much of the work needed for this is already done in the Coq mechanization
of Section~\ref{sectmechanize}.
The reason that this direction would be exciting is that the bounded-quantifier elimination
of Section~\ref{sect4} could be strengthened to full quantifier-elimination.

Finally, we are interested in embedding the hydra game of Kirby and Paris \cite{kirby}
(a short and very readable introduction is given by Bauer \cite{bauer})
into term interpretation.  A binary operator $+$ is (left and right)
\emph{self-distributive} if it satisfies
$a+(b+c)=(a+b)+(a+c)$ and $(a+b)+c=(a+c)+(b+c)$ (self-distributive operators were studied by
Frink \cite{frink} and more recently (left-sided only) by set theorists and knot theorists,
as surveyed by Dehornoy \cite{dehornoy}).
For such an operator (assuming also associativity),
\begin{align*}
{} &{}
\sum_{i_1=0}^{v_1} t_1 + \cdots + \sum_{i_k=0}^{v_k} t_k\\
&= \sum_{i_1=0}^{v_1} t_1 + \cdots + \sum_{i_{k-1}=0}^{v_{k-1}} t_{k-1} + t_k(i_k|0) + \cdots + t_k(i_k|v_k)\\
&=
\left(
\sum_{i_1=0}^{v_1} t_1 + \cdots + \sum_{i_{k-1}=0}^{v_{k-1}} t_{k-1} + t_k(i_k|0)\right)\\
{} &{\hphantom{=}}
+ \cdots +
\left(
\sum_{i_1=0}^{v_1} t_1 + \cdots + \sum_{i_{k-1}=0}^{v_{k-1}} t_{k-1} + t_k(i_k|v_k)\right),
\end{align*}
which bears a certain resemblance to the act of cutting a hydra's head and having many isomorphic copies of 
its subtree regrow.


\end{document}